# A Rank Revealing Randomized Singular Value Decomposition (R³SVD) Algorithm for Low-rank Matrix Approximations

Hao Ji, *Student Member, IEEE*, Wenjian Yu, *Member, IEEE*, and Yaohang Li, *Member, IEEE*

*Abstract*— In this paper, we present a Rank Revealing Randomized Singular Value Decomposition (R³SVD) algorithm to incrementally construct a low-rank approximation of a potentially large matrix while adaptively estimating the appropriate rank that can capture most of the actions of the matrix. Starting from a low-rank approximation with an initial guessed rank, R³SVD adopts an orthogonal Gaussian sampling approach to obtain the dominant subspace within the leftover space, which is used to add up to the existing low-rank approximation. Orthogonal Gaussian sampling is repeated until an appropriate low-rank approximation with satisfactory accuracy, measured by the overall energy percentage of the original matrix, is obtained. While being a fast algorithm, R³SVD is also a memory-aware algorithm where the computational process can be decomposed into a series of sampling tasks that use constant amount of memory. Numerical examples in image compression and matrix completion are used to demonstrate the effectiveness of R³SVD in low-rank approximation.

## I. INTRODUCTION

CONSIDERING an $m \times n$ matrix $A$ with rank $r$, the optimal $k$-rank ($k \leq r$) approximation $A_k$ of matrix $A$ yields minimum approximation error among all possible $m \times n$ matrices of rank $k$ [1], i.e.,

$$\|A - A_k\|_F^2 = \min_{rank(X)=k} \|A - X\|_F^2.$$

Within controllable approximation error, a good low-rank approximation of a large matrix can reduce storage requirement and accelerate matrix operations such as matrix-vector or matrix-matrix multiplications. If $A$ is a matrix representing data affinity in a large dataset, low-rank approximation algorithms can be used for dimensionality reduction or noise elimination. As a result, constructing appropriate low-rank approximations of large matrices plays a central role in many data analytic applications [2, 3, 4, 5, 6], such as principle component analysis, compressed sensing, data compression, signal processing, machine learning, and matrix completion.

The optimal $k$-rank approximation $A_k$ can be straightforwardly obtained by computing full Singular Value Decomposition (SVD) and then truncating it by selecting the top $k$ dominant singular values and their corresponding singular vectors such that

$$A_k = \sum_{i=1}^{k} \sigma_i u_i v_i^T,$$

where $k \leq r$, $\sigma_1$, $\sigma_2$, $\ldots$, $\sigma_k$ are the singular values of $A$ in non-increasing order, and $u_1, \cdots, u_k$ and $v_1, \cdots, v_k$ are the corresponding left and right singular vectors, respectively. Here, by tuning the value of $k$, the low-rank matrix approximation error measured by Frobenius norm can be controlled by

$$\|A - A_k\|_F^2 = \sum_{i=k+1}^{r} \sigma_i^2.$$

However, numerically computing the full SVD of a matrix when both $m$ and $n$ are large is often computationally costly as well as memory intensive. As an efficient alternatives, randomized algorithms to approximate SVD have attracted great interest recently and become competitive for computing rapid low-rank approximations of large matrices [3, 7, 8, 9]. Instead of passing over the large matrix in full SVD, the randomized SVD algorithms focus on efficiently sampling the important matrix elements. Many sampling strategies, including uniform column/row sampling (with or without replacement) [10, 11], diagonal sampling or column-norm sampling [12], sampling with $k$-means clustering [13], and Gaussian sampling [14], have been proposed. As a result, compared to full SVD, randomized SVD methods are memory efficient and can often obtain good low-rank approximations in a significantly faster way.

Nevertheless, most of these randomized SVD algorithms require the rank value $k$ to be given as an input parameter in advance. In fact, in many practical applications, $k$ is unknown beforehand but is of great importance to the accuracy of the solutions. In general, underestimating $k$ can introduce unacceptable large error in the low-rank approximation while

H. Ji and Y. Li are with the Department of Computer Science, Old Dominion University, Norfolk, VA 23529, USA (e-mail: hji@cs.odu.edu, yaohang@cs.odu.edu).

W. Yu is with Tsinghua National Laboratory for Information Science and Technology, the Department of Computer Science and Technology, Tsinghua University, Beijing 100084, China. (e-mail: yu-wj@tsinghua.edu.cn).

overestimating $k$ can lead to unnecessary computational and memory costs.

In this paper, we present a Rank-Revealing Randomized SVD (R³SVD) algorithm whose goal is to incrementally construct a low-rank approximation while estimating the appropriate $k$ value in an adaptive manner. The fundamental idea behind R³SVD is importance sampling – a new form of Gaussian sampling based on orthogonal projection is derived to obtain the dominant subspace orthogonal to the existing low-rank approximation, which is used to add up to existing low-rank approximation. Moreover, R³SVD is a memory-aware algorithm due to the fact that its computation can be tailored into subsequent tasks that can fit in constant amount of memory. We use several application examples, including image compression and matrix completion, to demonstrate the effectiveness of R³SVD.

The rest of the paper is organized as follows. We review the randomized SVD algorithm based on Gaussian Sampling for low-rank approximation in Section II. In Section III, we describe our R³SVD algorithm and justify its properties. Numerical examples are presented in Section IV. Finally, Section V summarizes our conclusions and future research directions.

## II. RANDOMIZED ALGORITHM WITH GAUSSIAN SAMPLING

Our R³SVD algorithm is based on the randomized SVD algorithm with Gaussian sampling proposed by Halko et al. [14, 15], although it can be straightforwardly extended to other randomized SVD algorithms with different matrix sampling strategies. In this paper, RSVD is referred to the randomized SVD algorithm with Gaussian sampling. The basic idea of RSVD is to use Gaussian vectors to construct a small condensed subspace from the range of $A$, whose dominant actions could be quickly estimated from this small subspace with relatively low computation cost while yielding high confidence. The procedure of RSVD is described in Algorithm 1.

---

**Algorithm 1: Randomized SVD (RSVD) Algorithm with Gaussian Sampling**

Input: $A \in \mathbb{R}^{m \times n}$, a target matrix rank $k \in \mathbb{N}$, and an oversampling parameter $p \in \mathbb{N}$ satisfying $k + p \leq min(m, n)$.
Output: Low rank approximation $U_L \in \mathbb{R}^{m \times k}$, $\Sigma_L \in \mathbb{R}^{k \times k}$, and $V_L \in \mathbb{R}^{n \times k}$

Construct an $n \times (k + p)$ Gaussian random matrix $\Omega$
$Y = A\Omega$
Compute an orthogonal basis $Q = qr(Y)$
$B = Q^T A$
$[U_B, \Sigma_B, V_B] = svd(B)$
Update $U_B = QU_B$
$U_L = U_B(:, 1:k)$, $\Sigma_L = \Sigma_B(1:k, 1:k)$, and $V_L = V_B(:, 1:k)$

---

Given a desired rank $k$ and an oversampling parameter $p$ (typically a small constant), RSVD constructs an $n \times (k + p)$ Gaussian random matrix block $\Omega$, whose elements are normally distributed. $\Omega$ condenses a large matrix $A$ into a "tall-and-skinny," dense block matrix $Y$ by $Y = A\Omega$. $Y$ captures the most important actions of $A$ and a basis $Q$ is derived by decomposing $Y$. $Q$ is designed to approximate the left singular vectors of $A$ by minimizing $||QQ^T A - A||_F^2$. Then, $Q$ is applied back to $A$ to obtain a "short-and-wide" block matrix $B = Q^T A$. Calculation of SVD on $B$ yields an approximate Singular Value Decomposition of $A$. The result $U_L \Sigma_L V_L^T$ forms a $k$-rank matrix approximation to $A$.

Compared to full SVD directly operating on the $m \times n$ matrix $A$, which is rather computational costly when both $m$ and $n$ are large, the major operations in RSVD are carried out on the block matrices instead. These block matrix operations include matrix-block matrix multiplications as well as QR and SVD decompositions on the block matrices. Specifically, matrix-block matrix multiplications take $O(2(k + p)T_{mult})$ floating-point operations, where $T_{mult}$ denotes the computational cost of a matrix-vector multiplication. For a large matrix $A$ where $m, n \gg k + p$, the computational cost of matrix-block matrix multiplications dominates those of block QR or SVD decomposition operations, which requires $O((k + p)^2 (m + n))$ floating operations. RSVD needs to store the intermediate matrices, such as $\Omega$, $Y$, $Q$, and $B$, and thus its space complexity is $O(2(m + n)(k + p))$. As a result, RSVD is usually more efficient than the full SVD algorithms in terms of computational and memory cost, but with a tradeoff of accuracy.

The desired rank $k$ is a required input parameter in the randomized SVD algorithms. However, in many practical applications, the value of $k$ is unknown beforehand and needs to be appropriately estimated. In the literature, two strategies have been proposed to estimate the appropriate value of $k$. One strategy is based on preprocessing, which intends to obtain an appropriate $k$ value before carrying out RSVD. For instance, Voronin and Martinsson [16] proposed two algorithms, Autorank I and Autorank II, to evaluate a basis $Q$ for a range space that captures the most actions of matrix $A$. Autorank I is based on overestimation by using a very large value $k$ at the beginning and then selecting dominant information from the resulting pool of singular values/vectors. Although Autorank I is often able to obtain good low-rank approximations, largely overestimated $k$ will result in significant computational cost increase, because the computational cost of decomposing the tall-and-skinny or short-and-wide block matrices in RSVD grows rapidly with $O(k^2)$ and is no longer negligible. At the same time, the memory requirement of Autorank I increases in the order of $k$, too. Instead of overestimating $k$, Autorank II gradually samples the range of $A$ guided by error $||QQ^T A - A||_F^2$ in order to obtain a good estimation of $k$. Similar to Autorank II, the Adaptive Randomized Range Finder algorithm [14] employs the incremental sampling approach with a probabilistic error estimator based on the relation between the rank $k$ with respect to the theoretical error bound to predict a reasonable basis $Q$ with a reasonable value of $k$. However, this theoretical error bound is loose and consequently $k$ is often largely overestimated, which will be shown in section 4. More recently, the Randomized Blocked algorithm [17], a block version of Randomized Range Finder algorithm, is developed to improve computational efficiency. Instead of using the probabilistic error estimator, the Randomized Blocked algorithm explicitly updates $A$ by removing the portion projected on $Q$ and terminates at a situation when $|A|$ becomes small enough. A practical issue of the Randomized Blocked algorithm is the "fill-in" problem – when the original matrix $A$ is sparse, explicitly updating $A$ will lead to emergence of a lot of non-zero elements which make $A$ become dense. Alternative to preprocessing, another strategy to obtain a good estimation of $k$ is to adaptively increase $k$ and evaluate its appropriateness while carrying out RSVD. A simple approach is restarting RSVD, which starts with a small guessed rank $k$ and then repeats RSVD computation with increasing $k$ until the low-rank approximation with desired accuracy is reached. This restarting approach can often result in a good low-rank approximation; however, the previous RSVD trials are only used to estimate $k$ and do not contribute to final low-rank approximation. The R³SVD algorithm presented in this paper is based on the second strategy. Unlike restarting RSVD, the low-rank approximation is built incrementally and therefore the previous RSVD trials are not wasted.

## III. RANK REVEALING RANDOMIZED SVD (R³SVD) ALGORITHM

To illustrate the R³SVD algorithm for a given matrix $A \in \mathbb{R}^{m \times n}$, we first define the energy $E(U_L \Sigma_L V_L^T)$ of a $k$-rank approximation $U_L \Sigma_L V_L^T$ as $E(U_L \Sigma_L V_L^T) = \|U_L \Sigma_L V_L^T\|_F^2$ while the overall energy of $A$ is $E(A) = \|A\|_F^2 = \sum_{i=1}^r \sigma_i^2$. Then, the overall energy percentage occupied by the $k$-rank approximation with respect to that of $A$ becomes

$$\varphi = \frac{\|U_L \Sigma_L V_L^T\|_F^2}{\|A\|_F^2} = \frac{\sum_{i=1}^k \sigma_{L_i}^2}{\sum_{i=1}^r \sigma_i^2}.$$

where $\sigma_{L_i}$ denotes the $i$th diagonal element of $\Sigma_L$ and $\sigma_i$'s are the actual singular values of $A$. Measuring the percentage of energy of a low-rank approximation with respect to a large matrix has been popularly used in a variety of applications for dimensionality reduction such as Principle Component Analysis (PCA) [4, 18], ISOMAP learning [19], Locally Linear Embedding (LLE) [20], and Linear Discriminant Analysis (LDA) [21]. According to the Eckart-Young-Mirsky theorem [22], for a fixed $k$ value, the optimal $k$-rank approximation has the overall energy percentage of $\sum_{i=1}^k \sigma_i^2 / \|A\|_F^2$ of $A$.

The adaptivity of R³SVD is achieved by estimating the overall energy percentage of the low-rank approximation obtained so far. The rationale of R³SVD is to build a low-rank approximation incrementally based on orthogonal Gaussian projection. Initially, a $t$-rank approximation is obtained, where $t$ is an initial guess of $k$ which can be justified according to the memory available. The energy percentage is estimated accordingly. If the energy percentage obtained so far does not satisfy the application requirement, a new $t$-rank approximation is calculated in the subspace orthogonal to the space of the previous low-rank approximation. Then, the new $t$-rank approximation will be added to the previous one to form a $2t$-rank approximation and its corresponding energy percentage is estimated. The above process is repeated until the incrementally built low-rank approximation has secured satisfactory percentage of energy from $A$.

Compared to RSVD, R³SVD incorporates three major changes including orthogonal Gaussian sampling, orthogonalization process, and stopping criteria based on energy percentage estimation. By putting all pieces together, R³SVD is described as follows.

---

**Algorithm 2: Rank Revealing Randomized SVD Algorithm**

---

**Input:** $A \in \mathbb{R}^{m \times n}$, sampling size $t \in \mathbb{N}$ per iteration, oversampling number $p \in \mathbb{N}$, power number $q \in \mathbb{N}$, maximum number of iterations $maxit \in \mathbb{N}$, and energy threshold $\tau \in \mathbb{R}$.
**Output:** Low rank approximation $U_L \in \mathbb{R}^{m \times k'}$, $\Sigma_L \in \mathbb{R}^{k' \times k'}$, $V_L \in \mathbb{R}^{n \times k'}$, and estimated rank $k'$

*// initialization*
Construct an $n \times (t + p)$ standard Gaussian matrix $\Omega$
$G_0 = \Omega$ and $V_L = \emptyset$, $U_L = \emptyset$, $\Sigma_L = \emptyset$
$k' = 0$
**for** $i = 0: maxit$
    $Y_i = AG_i$
    $Q_i = qr(Y_i, 0)$
    $B_i = Q_i^T A$
    $[U_{B_i}, \Sigma_{B_i}, V_{B_i}] = svd(B_i, 0)$
    $U_{B_i} = Q_i U_{B_i}$
    $V_{B_i} = qr(V_{B_i} - V_L(V_L^T V_{B_i}), 0)$    *// orthogonalization process*
    $U_L \leftarrow [U_L, U_{B_i}(:, 1:t)]$, $\Sigma_L \leftarrow \begin{bmatrix} \Sigma_L & 0 \\ 0 & \Sigma_{B_i}(1:t, 1:t) \end{bmatrix}$, $V_L \leftarrow [V_L, V_{B_i}(:, 1:t)]$
    **for** $j = 1: t$
        $k' = i \times t + j$
        $\tilde{\varphi}_{k'} = \frac{\sum_{i=1}^{k'} \sigma_{L_i}^2}{\|A\|_F^2}$    *// estimate energy percentage*
        **if** $\tilde{\varphi}_{k'} \geq \tau$, **then** stop;
    **end**
    $G_{i+1} = G_i - V_i(V_i^T G_i)$    *// update Gaussian matrix*
**end**
$[\Sigma_L, \text{Idx}] = \text{sort}(\Sigma_L, \text{'descend'});$    *// sort the approximate singular values*
$V_L = V_L (:, \text{Idx});$
$U_L = U_L (:, \text{Idx});$

---

### A. Orthogonal Gaussian Sampling

Suppose that $V_L$ is an $n \times t$ matrix composed of $t$ right singular vectors of a low-rank approximation $U_L \Sigma_L V_L$, which is supposed to capture most of the energy in $A$. Then, the range space, $ran(A^T)$, can be divided into two orthogonal spaces: space $ran(V_L)$ spanned by the columns in $V_L$ and its orthogonal complement $ran(V_L)^\perp$. Obviously, if $V_L$ consists of only partial dominant actions of $A$, the rest dominant information is left over in the space $ran(V_L)^\perp$.

R³SVD is designed to incrementally add up a low-rank approximation. Therefore, R³SVD needs to sample the space $ran(V_L)^\perp$ orthogonal to $V_L$ to extract the left-over dominant information of $A$. Here, we construct a sampling matrix $G$ such as

$$G = (I - P_V)\Omega$$

where $P_V = V_L V_L^T$ is an orthogonal projection onto the space $ran(V_L)$, $\Omega$ is a standard Gaussian matrix, and $I$ is the identity matrix. Theorem 1 shows that $G$ is a Gaussian matrix orthogonal to $ran(V_L)$.

**Theorem 1.** Assuming that $V_L$ is an $n \times t$ non-empty matrix with orthonormal columns, then

1) $G$ is orthogonal to $V_L$; and

2) elements in $G$ are normally distributed.

*Proof.* 1) Since $V_L$ is an $n \times t$ matrix with orthonormal columns, $P_V$ can be derived as $P_V = V_L V_L^T$. Obviously, $V_L^T(I - P_V) = V_L^T - V_L^T V_L V_L^T = 0$ holds.

2) As $I - P_V$ is the orthogonal projection onto $ran(V_L)^\perp$, which is the orthogonal complement of space $ran(V_L)$, we can denote an $n \times (n-t)$ matrix $\tilde{V} = (\tilde{v}_{ij})$ as a basis of the space $ran(V_L)^\perp$ and then $I - P_V = \tilde{V}\tilde{V}^T$. Then, each element $g_{ij}$ in $G$ can be expressed as

$$g_{ij} = \sum_{s=1}^{n}\left(\sum_{h=1}^{n-t}\tilde{v}_{ih}\tilde{v}_{sh}\right)\omega_{sj}$$

where $\omega_{sj}$ denotes an element of $\Omega$ in row $s$ of column $j$. Since element $\omega_{sj}$'s are independent standard normal distributed variables, the characteristic function $\Phi_{g_{ij}}(x)$ can be obtained as

$$\Phi_{g_{ij}}(x) = \Phi_{\sum_{s=1}^{n}(\sum_{h=1}^{n-t}\tilde{v}_{ih}\tilde{v}_{sh})\omega_{sj}}(x)$$
$$= \prod_{s=1}^{n}\prod_{h=1}^{n-t}e^{-\frac{1}{2}(\tilde{v}_{ih}\tilde{v}_{sh}x)^2}$$
$$= \prod_{h=1}^{n-t}\prod_{s=1}^{n}e^{-\frac{1}{2}(\tilde{v}_{ih}\tilde{v}_{sh}x)^2}$$
$$= \prod_{h=1}^{n-t}e^{-\frac{1}{2}\tilde{v}_{ih}^2(\sum_{s=1}^{n}\tilde{v}_{sh}^2)x^2}.$$

As the columns of $\tilde{V}$ are orthonormal such that $\left(\sum_{s=1}^{n}\tilde{v}_{sh}^2\right) = 1$, we have

$$\Phi_{g_{ij}}(x) = e^{-\frac{1}{2}\sum_{h=1}^{n-t}\tilde{v}_{ih}^2 x^2}.$$

Since the characteristic function uniquely determines the probability distribution of a random variable [23], it suffices to show that $g_{ij}$ is normally distributed with expected value zero and variance $\tilde{\sigma}_{ij}^2 = \sum_{h=1}^{n-t}\tilde{v}_{ih}^2$, i.e., $g_{ij} \sim N(0, \tilde{\sigma}_{ij}^2)$.

To avoid resampling of the space $ran(V_L)$, the product of matrix $AG$ in R³SVD focuses on revealing the dominant information from the space $ran(V_L)^\perp$ orthogonal to $ran(V_L)$. Since the number of dominant singular values is unknown in advance, R³SVD generates a series of Gaussian matrices $G_1, G_2, ...$ to iteratively explore the orthogonal subspace of the obtained low-rank approximation until a satisfactory low rank approximation is obtained.

In the situation when the singular spectrum of matrix $A$ decays slowly, Gaussian matrix can be applied to a power multiplication of matrix $A$ to improve the approximation accuracy. This is referred to as the power scheme suggested in [14]. In R³SVD, the power scheme with Gaussian matrices $G_i$ is applied to $(A(I - P_V)A^T)^q A$ where $q$ is a power scalar. This power scheme is able to refine the sampled space $Y_i$ orthogonal to space $ran(V_L)$. Algorithm 3 shows the R³SVD algorithm with power scheme.

---

**Algorithm 3: Rank Revealing Randomized SVD Algorithm with Power Scheme**

---

**Input:** $A \in \mathbb{R}^{m \times n}$, sampling size $t \in \mathbb{N}$ per iteration, oversampling number $p \in \mathbb{N}$, power number $q \in \mathbb{N}$, maximum number of iterations $maxit \in \mathbb{N}$, and energy threshold $\tau \in \mathbb{R}$.
**Output:** Low rank approximation $U_L \in \mathbb{R}^{m \times k'}$, $\Sigma_L \in \mathbb{R}^{k' \times k'}$, $V_L \in \mathbb{R}^{n \times k'}$, and estimated rank $k'$

*// initialization*
Construct an $n \times (t + p)$ standard Gaussian matrix $\Omega$
$G_0 = \Omega$ and $V_L = \emptyset$, $U_L = \emptyset$, $\Sigma_L = \emptyset$
$k' = 0$
**for** $i = 0 : maxit$
    $Y_i = AG_i$
    $Q_i = qr(Y_i, 0)$

    **for** $j = 1 : q$                *// the power scheme*
        $Y_i = A^T Q_i$
        $Y_i = Y_i - V_L(V_L^T Y_i)$
        $Q_i = qr(Y_i, 0)$
        $Y_i = AQ_i$
        $Y_i = Y_i - V_L(V_L^T Y_i)$
        $Q_i = qr(Y_i, 0)$
    **end**

    $B_i = Q_i^T A$
    $[U_{B_i}, \Sigma_{B_i}, V_{B_i}] = svd(B_i, 0)$
    $U_{B_i} = Q_i U_{B_i}$
    $V_{B_i} = qr(V_{B_i} - V_L(V_L^T V_{B_i}), 0)$        *// orthogonalization process*
    $U_L \leftarrow [U_L, U_{B_i}(:, 1:t)]$, $\Sigma_L \leftarrow \begin{bmatrix} \Sigma_L & 0 \\ 0 & \Sigma_{B_i}(1:t, 1:t) \end{bmatrix}$, $V_L \leftarrow [V_L, V_{B_i}(:, 1:t)]$
    **for** $j = 1 : t$
        $k' = i \times t + j$
        $\tilde{\varphi}_{k'} = \frac{\sum_{i=1}^{k'}\sigma_i'^2}{\|A\|_F^2}$        *// estimate energy percentage*
        **if** $\tilde{\varphi}_{k'} \geq \tau$, **then** stop;
    **end**
    $G_{i+1} = G_i - V_i(V_i^T G_i)$    *// update Gaussian matrix*

**end**
$[\Sigma_L, \text{Idx}] = \text{sort}(\Sigma_L, '\text{descend}');$     *// sort the approximate singular values*
$V_L = V_L(:, \text{Idx});$
$U_L = U_L(:, \text{Idx});$

### B. Orthogonalization Process

Let $V_L = [V_1, V_2, \ldots V_i]$ denote a matrix containing the approximate right singular vectors obtained in R³SVD after the $i$th iteration step. Then, the singular vectors in $V_{i+1}$ must be orthogonal to $V_L$. However, the inherent numerical errors may cause loss of orthogonality between $V_{i+1}$ and $V_L$.

To ensure the orthogonality property, we generate $V_{i+1}$ by employing an orthogonalization process to remove the components of $V_{B_i}$ that are not orthogonal to the previous right singular vectors in $V_L$ such that

$$V_{i+1} = qr(V_{B_i} - V_L(V_L^T V_{B_i}), 0).$$

Proposition 2 indicates that the resulting matrix $V_{i+1}$ generated at the $(i+1)$th iteration step in R³SVD is orthogonal to $V_L$.

**Proposition 2.** $V_L^T V_{i+1} = 0$ holds.
*Proof.* Given $Z_{i+1} = V_{B_i} - V_L(V_L^T V_{B_i})$, we can get $V_L^T Z_{i+1} = V_L^T(V_{B_i} - V_L(V_L^T V_{B_i})) = 0$.
Since $V_{i+1}$ is a basis of $ran(Z_{i+1})$, $V_L^T V_{i+1} = 0$ holds.

Based on Proposition 2, the orthogonality property of the resulting left singular vectors $U_L$ is proved in Proposition 3.

**Proposition 3.** $U_L^T U_{i+1} = 0$ holds.
*Proof.* Denoting the QR decomposition of $Y_{i+1}$ by $Y_{i+1} = Q_{i+1} R_{i+1}$. We can have
$$\begin{aligned}
Q_j^T Q_{i+1} &= Q_j^T Y_{i+1} R_{i+1}^{-1} \\
&= Q_j^T A G_{i+1} R_{i+1}^{-1} \\
&= Q_j^T A(I - P_{V_L}) \Omega R_{i+1}^{-1} \\
&= B_j(I - P_{V_L}) \Omega R_{i+1}^{-1} \\
&= U_{B_j} \Sigma_{B_j} V_{B_j}^T (I - P_{V_L}) \Omega R_{i+1}^{-1}
\end{aligned}$$
where $V_L = [V_1, V_2, \ldots V_i]$.
Denote $V^- = [V_1, V_2, \ldots V_{j-1}]$ and $V^+ = [V_{j+1}, \ldots, V_i]$ for $j \leq i$. According to Proposition 2 that the columns in $V_L$ are orthogonal to each other, $(I - P_{V_L})$ can be expressed as
$$(I - P_{V_L}) = (I - P_{V^-})(I - P_{V_j})(I - P_{V^+}).$$
Since $V_j = qr(V_{B_j} - V^-(V^{-T} V_{B_j}), 0)$, it follows that
$$(I - P_{V^-}) V_{B_j} = V_j R_j.$$
Therefore,
$$\begin{aligned}
V_{B_j}^T(I - P_{V_L}) &= V_{B_j}^T(I - P_{V^-})(I - P_{V_j})(I - P_{V^+}) \\
&= R_j^T V_j^T(I - P_{V_j})(I - P_{V^+}) \\
&= 0
\end{aligned}$$
and thus $Q_j^T Q_{i+1} = 0$, for $j \leq i$. In conclusion,
$$U_L^T U_{i+1} = [U_1, U_2, \ldots U_i]^T U_{i+1}$$

$$\begin{aligned}
&= \begin{bmatrix} U_{B_1}^T Q_1^T Q_{i+1} U_{B_{i+1}} \\ U_{B_2}^T Q_2^T Q_{i+1} U_{B_{i+1}} \\ \vdots \\ U_{B_{i+1}}^T Q_i^T Q_{i+1} U_{B_{i+1}} \end{bmatrix} \\
&= 0.
\end{aligned}$$

The orthogonalization process requires $O((2ti+1)(t+p)n)$ operations to ensure the orthogonality properties of singular vectors obtained in the previous iterations. Moreover, by taking advantage of the orthogonality between $V_{i+1}$ and $V_L$, the next Gaussian matrix $G_{i+1}$ can be fast generated using the following short recursive formula,

$$G_{i+1} = \left(I - \sum_{j=1}^{i} P_{V_j}\right)\Omega = G_i - P_{V_i} G_i,$$

where $P_{V_j}$ is the orthogonal projection onto the space spanned by $V_j$, such that $P_{V_j} = V_j V_j^T$, and $\Omega$ is a standard Gaussian matrix. Since $G_{i+1}$ is generated directly from $G_i$, the orthogonal Gaussian sampling takes only $O((2t+1)(t+p)n)$ operations.

### C. Energy Estimation and Stopping Criteria

The incremental low-rank approximation buildup process in R³SVD will be terminated when sufficient percentage of energy of $A$ is secured. The energy percentage threshold $\tau$ is typically specified by the applications, which often ranges from 80% to 99%.

Let $U_L = [U_1, U_2, \ldots U_i]$ denote a matrix of the approximate left singular vectors. The actual energy percentage of the low-rank approximation obtained at the $i$th iteration step can be evaluated based on

$$\varphi = \frac{\|U_L U_L^T A\|_F^2}{\|A\|_F^2}$$

However, to avoid costly calculation of $\|U_L U_L^T A\|_F^2$ at each iteration, in this paper, we adopt the following measure $\tilde{\varphi}_{k'}$ to quickly estimate the energy percentage of the obtained low-rank approximation. Here

$$\tilde{\varphi}_{k'} = \frac{\sum_{i=1}^{k'} \sigma_i'^2}{\|A\|_F^2}.$$

where $\sigma_i'$ denotes the $i$th approximate singular value in R³SVD. Proposition 4 shows that the estimated energy $\tilde{\varphi}_{k'}$ is equivalent to the actual energy $\varphi$. Therefore, it guarantees that the low-rank approximation obtained by R³SVD satisfies the accuracy requirement of the applications.

**Proposition 4.** $\tilde{\varphi}_{k'} = \varphi$
Proof. Since the columns in $U_L$ are orthogonal, we have
$$\begin{aligned}
\|U_L U_L^T A\|_F^2 &= \sum_{j=1}^{i} \|U_j U_j^T A\|_F^2 \\
&= \sum_{j=1}^{i} \|U_j^T A\|_F^2
\end{aligned}$$

$$= \sum_{j=1}^{i} \left\| U_{B_j}^{\ T} Q_j^{\ T} A \right\|_F^2$$

$$= \sum_{j=1}^{i} \left\| U_{B_j}^{\ T} U_{B_j} \Sigma_{B_j} V_{B_j}^{\ T} \right\|_F^2$$

$$= \sum_{j=1}^{i} \left\| \Sigma_{B_j} \right\|_F^2$$

$$= \sum_{i=1}^{k'} \sigma_i'^2$$

Hence,

$$\tilde{\varphi}_{k'} = \frac{\sum_{i=1}^{k'} \sigma_i'^2}{\|A\|_F^2} = \varphi' = \frac{\|U_L U_L^T A\|_F^2}{\|A\|_F^2}.$$

It is important to note the approximate singular values $\sigma_i'$'s are available during the calculation of SVD on $B_i$, where $B_i = Q_i^T A$. Therefore, the energy percentage can be evaluated at (almost) no cost.

### D. Complexity Analysis

As discussed above, at each iteration, R$^3$SVD carries out orthogonal Gaussian sampling to compute a new $t$-rank approximation of the leftover subspace orthogonal to the low-rank approximation obtained so far. Suppose that R$^3$SVD uses $s$ iterations to achieve a satisfactory low rank approximation with $k' \approx ts$ as the result rank and assume that the computational cost of matrix-block matrix multiplications dominates those of QR and SVD decompositions on the block matrices. The computational cost of R$^3$SVD is

$$O(2(k' + sp)T_{mult}).$$

In the case that matrix $A$ is sparse and both $m$ and $n$ are large, we can obtain the time complexity with simpler terms. In particular, as $T_{mult} \approx cm$, where $c$ is sparsity ratio, the time complexity of R$^3$SVD can be expressed as $O(k'^2 min(m, n))$.

In additional to the storage of matrix $A$, the major computations of R$^3$SVD are carried out on a series of block matrices with $(t + p)$ columns or rows. Therefore, where $t < k$, R$^3$SVD takes a constant space complexity of $O(2(m + n)(t + p))$, which is lower than that of RSVD, $O(2(m + n)(k + p))$.

## IV. NUMERICAL RESULTS

In this section, we use several numerical examples to demonstrate the effectiveness of R$^3$SVD for low-rank approximation in image compression and matrix completion.

### A. Comparisons with RSVD

We compare the performance of R$^3$SVD, full SVD, Autorank II, restarting RSVD, Adaptive Randomized Range Finder algorithm, and Randomized Blocked algorithm in constructing low rank approximations to compress a $7671 \times 7680$ NASA synthesis image chosen from the

Mars Exploration Rover mission [24]. The energy percentage threshold $\tau$ is set to 99%.

Both R$^3$SVD and RSVD start with an initial guess $t = 15$ of the target rank and $p = 5$ extra oversampling vectors. The power scheme is not applied, such that $q = 0$. In restarting RSVD, as the approximate singular values are available during each RSVD trial, the energy estimation introduced in Section III is used. If the guessed rank turns out to be insufficient to obtain a low rank approximation with satisfactory accuracy, the restarting approach repeats the RSVD computation with a gradually increasing rank $\Delta t$=15. Table I compares the computational performance of full SVD, Autorank II algorithm, Adaptive Randomized Range Finder algorithm, Randomized Blocked Algorithm, R$^3$SVD, and restarting RSVD in terms of rank, computational time, maximum memory usage, and energy percentage of the obtained low rank approximation. The optimal low-rank approximation (rank 46) to obtain 99% energy of the original matrix can be obtained by carrying out full SVD, which takes over 760 seconds on a Dell Precision-M6500 laptop (Intel CoreTM i5CPU, 2.67GHz, 4GBRAM). Restarting RSVD reduces the computational time to 13.77 seconds with a low-rank approximation of rank 79. Compared to restarting RSVD, R$^3$SVD further reduces both the computational time to 4.54 (32.97%) and rank to 62 (78.48%). This is because R$^3$SVD carries out important sampling based on the approximate right singular vectors in $V_L$, which is computed by multiplying $A$ twice per iteration. The power scheme allows more precise estimation of the dominant actions than a single iteration of $A$ multiplication in restarting RSVD. It is also important to notice that the algorithms based on the strategy of estimating $k$ before RSVD, including Autorank II, Adaptive Randomized Range Finder, and Randomized Blocked Algorithm require more computational time as well as the memory than restarting RSVD and R$^3$SVD. The Adaptive Randomized Range Finder uses a probabilistic error estimator, which leads to a highly overestimated rank (641).

TABLE I
PERFORMANCE COMPARISON OF R$^3$SVD, FULL SVD, AUTORANK II, RESTARTING RSVD, ADAPTIVE RANDOMIZED RANGE FINDER, AND THE RANDOMIZED BLOCKED ALGORITHM

| | Rank | Computational Time (second) | Maximum Memory Usage (bytes) | Energy Percentage Achieved |
|---|---|---|---|---|
| Full SVD | 46 | 760.55 | $1.41 \times 10^9$ | 99.024% |
| Autorank II Algorithm [16] | 105 | 32.66 | $2.03 \times 10^7$ | 99.184% |
| Adaptive Randomized Range Finder Algorithm [14] | 641 | 55.65 | $1.19 \times 10^8$ | 99.999% |
| Randomized Blocked Algorithm [17] | 105 | 20.90 | $4.80 \times 10^8$ | 99.184% |
| Restarting RSVD | 79 | 13.77 | $1.60 \times 10^7$ | 99.000% |
| R$^3$SVD | 62 | 4.54 | $4.91 \times 10^6$ | 99.006% |

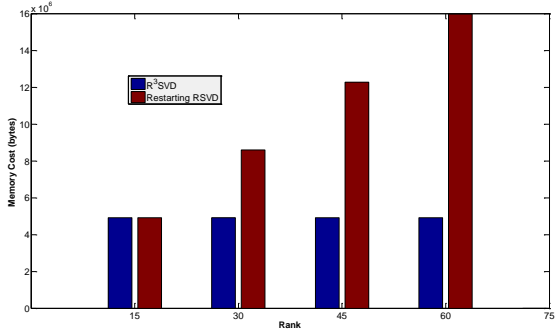

Fig. 1. Memory usage in R³SVD and restarting RSVD

Another important advantage of R³SVD is that R³SVD maintains constant memory usage in the computational process. Fig. 1 shows the memory usages in R³SVD and restarting RSVD as the guessed rank gradually increases. One can find that for a larger guessed rank, restarting RSVD requires more memory because of decomposing block matrices with more columns or rows. In contrast, the decomposition operations in R³SVD are carried out on block matrices with fixed $(t + p)$ number of columns or rows. As a result, the memory usage does not increase as the guessed rank increases during R³SVD. As shown in Table I, the memory usage in R³SVD is significantly less than those of the other algorithms.

Fig. 2 presents the compressed images in R³SVD, where Fig. 2(a) is the original image and Figs. 2(b) to 2(d) illustrate the adaptive compressed images with increasing ranks. With the resulting 55-rank low-rank approximation, a compressed image with 99.18% energy of the original image is obtained.

One advantage of R³SVD is that the computational process can be tailored into a series of sampling tasks that can fit into the available memory in a computer via adjusting the sampling size parameter $t$. Fig. 3 compares the memory usage in R³SVD with $t = 20$, 15, 10, and 5. One can find that a smaller sample size in R³SVD yields proportionally less consumption of memory but without significantly affecting the rank in the obtained low-rank approximation. The resulting ranks are 63, 62, 61, and 60, respectively. Therefore, calculating sampling size parameter $t$ according to the available memory in a computer can lead to the best computational performance of R³SVD.

### B. Application in Matrix Completion

R³SVD can be effectively applied to applications of matrix completion, whose goal is to recover the missing (unknown) entries of an incomplete matrix [6, 25, 26, 27]. Matrix completion algorithms have been widely used in many applications, including machine learning [28, 29], computer vision [30], and image/video processing [31]. Low rank matrix approximation is a core component in many matrix completion algorithms. The computational efficiency of constructing high-quality low rank approximation is essential

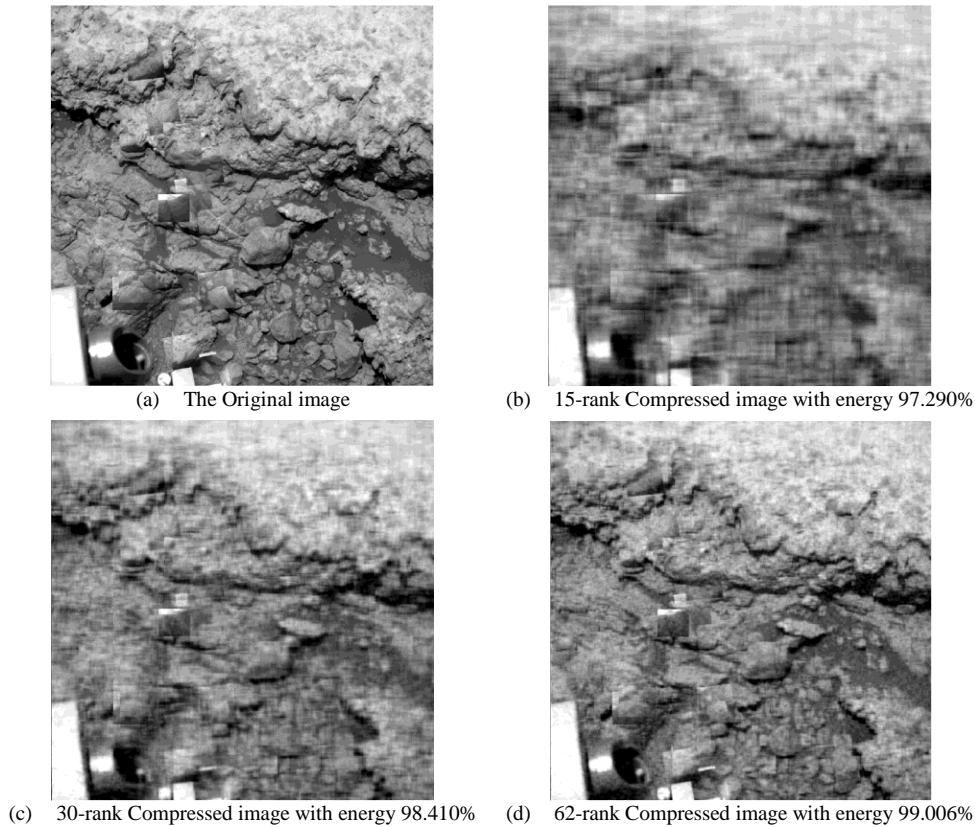

(a)    The Original image

(b)    15-rank Compressed image with energy 97.290%

(c)    30-rank Compressed image with energy 98.410%

(d)    62-rank Compressed image with energy 99.006%

Fig. 2. The original image and the compressed images with increasing ranks in R³SVD

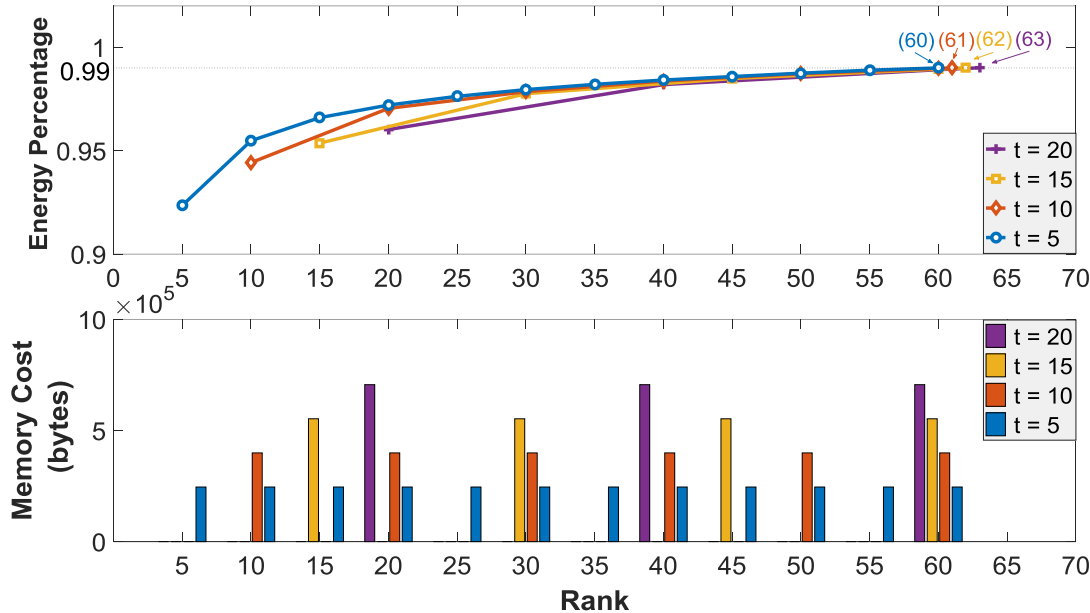

Fig. 3. The energy percentage of the obtained low rank approximations (upper) and the required memory space (lower) in R³SVD with $t = 20$, 15, 10, and 5 and the oversampling parameter $p = 5$.

to the performance of these matrix completion algorithms.

We modify the Singular Value Thresholding (SVT) algorithm [6] by replacing the underlying Lanczos algorithm with our R³SVD algorithm to compute dominant singular values and vectors at each SVT iteration. Fig. 4 shows a $1024 \times 1024$ aerial image chosen from the USC-SIPI Image Database [32] as well as 10% of the pixels uniformly sampled from the image (the background is set to grey to highlight these samples). As shown in Table II, the modified SVT algorithm obtains the completed image with similar recovery error and rank as that of the original SVT. In comparison, replacing Lanczos algorithm with R³SVD significantly reduces the overall computational time in SVT. This is due to the fact that the Lanczos bidiagonalization algorithm with partial reorthogonalization used in original SVT has computational complexity of $O(\min(m, n)^2 k)$ [33, 34] while R³SVD offers a faster way with computational complexity of $O(\min(m, n)k^2)$ in contrast. As a result, the modified SVT method using R³SVD achieves about 1.69 times speedup over the original SVT method using Lanczos algorithm.

## V. CONCLUSIONS

We present an R³SVD algorithm based on orthogonal sampling to gradually build up a low rank approximation of a given matrix to satisfy application specific accuracy. A random matrix based on the orthogonal complement operator is derived to enable R³SVD to concentrate on sampling the orthogonal subspace of the existing low-rank approximation. Compared to the algorithms based on preprocessing strategy by estimating the appropriate rank $k$ before RSVD, R³SVD is more efficient in terms of both computation time and memory while providing a better rank estimation. Moreover, as a memory-aware algorithm, R³SVD is particular favorable for many real-life applications running in limited computer memory. The effectiveness of R³SVD has been demonstrated in numerical applications including image compression and matrix completion.

The importance sampling approach proposed in this paper can also be used for other randomized algorithms by sampling the most important subspaces toward the solutions. The R³SVD algorithm described in this paper is based on Gaussian sampling [14, 15], which can also be extended to other randomized SVD strategies such as column or row sampling [10, 11].

## ACKNOWLEDGMENT

Yaohang Li acknowledges support from National Science Foundation grant number 1066471. Hao Ji acknowledges support from Old Dominion University Modeling and Simulation Fellowship. We acknowledge Shidong Jiang for stimulating discussions.

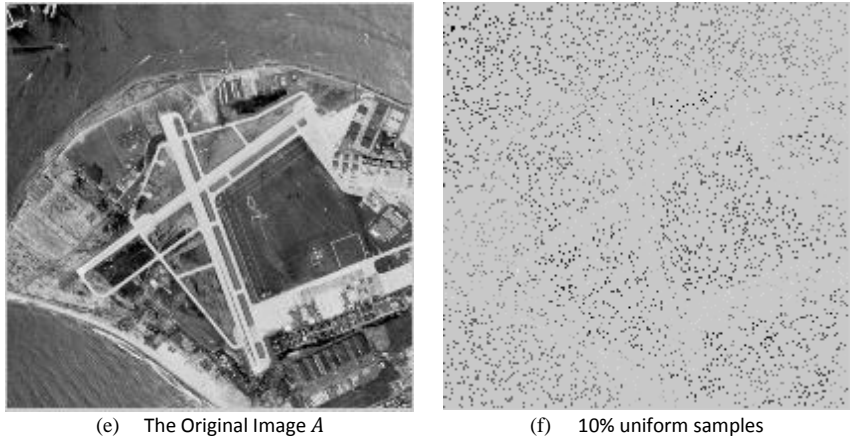

(e)  The Original Image $A$          (f)  10% uniform samples

Fig. 4. The original image and the sample image

TABLE II
THE COMPLETED IMAGES USING THE ORIGINAL SVT ALGORITHM AND THE MODIFIED SVT ALGORITHM USING R³SVD

| | Completed Image $X$ | Elapsed Time (seconds) | Rank | $\dfrac{\|A-X\|_F^2}{\|A\|_F^2}$ |
|---|---|---|---|---|
| Original SVT | 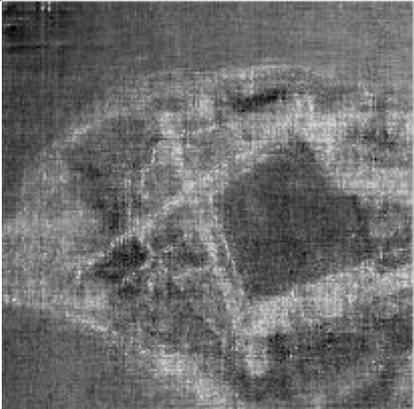 | 3457.122 | 190 | $6.772 \times 10^{-2}$ |
| Modified SVT using R³SVD | 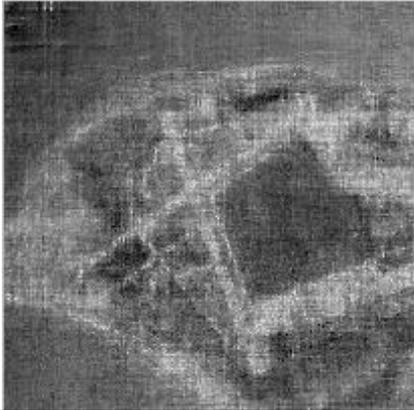 | 2045.295 | 189 | $6.772 \times 10^{-2}$ |